\documentclass[12pt]{amsart}
\usepackage{latexsym, amssymb}
\def\Tor{{\rm Tor}}
\def\proof{\medskip\noindent{\sc Proof. }}
\def\EOP{\hfill$\Box$}

\def\real{{\mathbb R}}

\newtheorem{theorem}{Theorem}
\newtheorem{lemma}[theorem]{Lemma}

\newtheorem{definition}[theorem]{Definition}

\newtheorem{question}[theorem]{Question}
\newtheorem{example}[theorem]{Example}

\newtheorem{remark}[theorem]{Remark}

\begin{document}

\title{Coloring complexes and arrangements}

\author{Patricia Hersh}
\address{Department of Mathematics, Indiana University, Rawles Hall, Bloomington, IN 47405}
\email{phersh@indiana.edu}
\author{Ed Swartz}
\address{Department of Mathematics, Cornell University, Ithaca, NY, 14853}
\email{ebs@math.cornell.edu}

\thanks{The first author was supported by NSF grant
DMS-0500638.  The second author was supported by NSF grant DMS-0245623}

\begin{abstract}

Steingrimsson's coloring complex and Jonsson's unipolar complex are interpreted in terms of hyperplane arrangements.

This viewpoint leads to short proofs 
that all coloring complexes  and a large class of unipolar complexes 
have convex ear decompositions.  These convex ear decompositions 
impose strong new restrictions on the chromatic polynomials of all finite graphs.  Similar results are obtained for characteristic polynomials of submatroids of type $\mathcal{B}_n$ arrangements.

\end{abstract}

\maketitle

\section{Introduction}

Since its introduction by Birkhoff  almost a century ago \cite{Bi}, the chromatic polynomial has been the object of intense study.  Nonetheless, a satisfactory answer to Wilf's question, ``What polynomials are chromatic?" \cite{Wi} remains elusive.  In \cite{Jo}, Jonsson proved that Steingrimsson's coloring complex is Cohen-Macaulay, and thereby established new restrictions on such polynomials.  Our main result is that the coloring complex has a convex ear decomposition, which implies that the chromatic polynomials of all finite graphs satisfy much stronger inequalities than those provided by \cite[Theorem 1.4]{Jo}.   

We also apply our methods to   Jonsson's unipolar complex and to characteristic polynomials of submatroids of type $\mathcal{B}_n$ arrangements.    On the other hand, we give examples indicating that these results cannot be extended to the characteristic polynomials of all matroids or even to large classes 
that seem to be particularly natural candidates.

The coloring complex $\Delta_G$ of a graph was introduced in [Ste] and was 
proven to be constructible, hence Cohen-Macaulay,
in [Jo].   The $(r-1)$-dimensional faces of the coloring complex are ordered lists
$T_1 |T_2|T_3| \cdots | T_r$ of nonempty disjoint sets of vertices with the property that at 
least one $T_i$ includes a pair of vertices that comprise an edge of $G$ and 
$\cup_{1\le i\le r}  T_i \ne V(G)$.  Steingrimsson showed that the $h$-polynomial of the 
double cone of
the coloring complex is related to the chromatic polynomial by the following formula.
\begin{equation} 
(1-t)^n \sum^\infty_{j=0} [(j+1)^n - P_G(j+1)] t^j = h_0 + h_1 t + \dots + h_n t^n.
\end{equation}
This expression  allows any new constraints on the $h$-vector of the coloring complex to be 
translated into new constraints on chromatic polynomials of all finite graphs.
Steingrimsson proved this formula by 
a Hilbert series calculation, so next we describe the rings involved.

Following \cite{Ste}, let $G$ be a graph with vertex set $V=[n]$. Set 
$A = k[x_S | S\subseteq [n]]$, $I = \langle x_Sx_T | S\not\subseteq T 
\mbox{ and } T \not\subseteq S \rangle $, and let
$R = A/I.$  By definition, $R= k[\Delta (B_n)], $ the Stanley-Reisner ring of the order complex of the Boolean algebra $B_n.$ Let $K_G$ be the ideal in $R$ generated by 
monomials $x_{S_1}x_{S_2}\cdots x_{S_r}$ such that for each $i\ge 1$ we have 
that $S_i \setminus S_{i-1}$ does not
include any pairs $\{i_1,i_2\}$ in $E(G),$ the edge set of $G$.  By convention, $S_0 = 
\emptyset $ so that $S_1 \setminus S_0 = S_1$ must be a disconnected set of vertices.   $K_G$ is often called the coloring ideal
of $G$.  It turns out that  $R/K_G$ is the 
Stanley-Reisner ring of the double cone of $\Delta_G$.  

In [Br], Brenti asked whether there exists, for an arbitrary graph $G$, 
a standard graded algebra whose Hilbert polynomial is the chromatic 
polynomial of $G$.   In general it is not possible for the Hilbert function of a standard graded algebra to agree identically with the values of the chromatic polynomial of a graph since the  latter is zero below the graph's chromatic number.  However, Steingrimsson showed that $K_G$ is an ideal whose Hilbert function agrees (up to a shift of one) with the values of the chromatic polynomial \cite{Ste}, and 
thereby obtained the above formula as a corollary.  In \cite{Ste}, he also 
attributes to G. Almkvist an earlier, nonconstructive affirmative answer to Brenti's question.

Steingrimsson's idea was to give a correspondence between the monomials in $K_G$ of 
degree $r$ and the proper $r+1$ colorings of $G$ as follows: the monomial 
$(x_{S_1})^{d_1}\cdots (x_{S_l})^{d_l}$ corresponds to the coloring in which the 
vertices in $S_1$ are colored $1$, the vertices in $S_2\setminus S_1$ are colored
$d_1 + 1$, the vertices in $S_3\setminus S_2$ are colored $d_1+d_2 +1$, etc.  Note that 
$S_1 = \emptyset $ if no vertices are colored 1.  We then have
$r = \sum d_i $, in other words, the degree of the monomial.  

In addition to proving that coloring complexes are constructible in [Jo], Jonsson also introduced the unipolar complex, proved it to be constructible, and determined its homotopy type.
By examining these complexes from the viewpoint of hyperplane arrangements we will prove that the coloring complex has a convex ear decomposition and that if the graph contains a vertex of degree $n-1$, then the unipolar complex also has a convex ear decomposition.  .  
From these results, we obtain new restrictions on the chromatic polynomials of all finite graphs in 
Section \ref{enum-consequences}.   See Section \ref{shell+ce} for the definition of convex ear decomposition. Applying this idea to subarrangements of type $\mathcal{B}_n$ arrangements leads to restrictions on their characteristic polynomials.

We assume the reader is familiar with Stanley-Reisner rings and $h$-vectors of finite simplicial complexes as presented in \cite{Sta}.  In Section \ref{matroids} we assume the reader is familiar with the characteristic polynomial of a matroid and its connection to the chromatic polynomial of a graph.  See, for instance, \cite[Section 6.3]{BO}

\section{An arrangements interpretation for the coloring complex}
Given a graph
$G$ with $n$ vertices, 
let $A_G$ be the real hyperplane arrangement generated by the hyperplanes
of the form $x_i = x_j$ for each edge $\{i,j\}$ present in $E(G)$.  When $G$ is $K_n,$ the complete graph on $n$ vertices, $A_{K_n}$ is usually called the type A braid arrangement.  In this case the intersection of all the hyperplanes is the line $x_1=x_2= \cdots = x_n.$  Let $H$ be the hyperplane $\{(a_1,\dots,a_n) \in \real^n: \sum a_i = 0\}.$   Then $A_{K_n} \cap H$ induces a simplicial cell decomposition on $S^{n-2},$ the unit sphere of $H$. The faces of the complex correspond to ordered  partitions  $S_1 | S_2 | \cdots |S_{r-1}| S_r,\  r \ge 2,$ of $[n].$  A point $(a_1,a_2,\dots ,a_n)$ is in the cell in which $S_1$ consists of those coordinates 
which are all equal 
to each other and are smaller than all other coordinates, and 
where $S_i$ is defined inductively to consist of all coordinates that are
 all 
equal to each other and are smaller than all other elements of
$\{ a_1,\dots ,a_n \} \setminus (S_1\cup\cdots \cup S_{i-1} )$. The top dimensional faces have dimension $n-2$ and correspond to partitions with $|S_i|=1$ for all $i.$   Identifying ordered partitions $S_1 | S_2 | \cdots | S_{r-1}| S_r$ of $[n]$ with ordered partitions $S_1 | S_2 | \cdots | S_{r-1}$ of proper subsets of $[n],$ the above discussion makes it clear that $\Delta_{K_n}$ is simplicially isomorphic to the codimension one skeleton of $S^{n-2} \cap A_{K_n}.$  In addition, we can see from its definition, that
$\Delta_G$ is isomorphic as a simplicial complex to the restriction of $A_{K_n} $ to
$(S^{n-2} \cap A_G).$  
The above discussion is essentially a special case of an idea appearing in [HRW].  We sum up the above with the following theorem.

\begin{theorem} \label{arrangement}
The coloring complex of $G$ is isomorphic as a
simplicial complex to the restriction 
of $A_{K_n} \cap S^{n-2}$ to the arrangement $A_G$.
\end{theorem}

One consequence is a new, short proof of the following result (also see
Theorem 4.2 of [HRW] for a generalization of this result). 

\begin{theorem}[Jonsson]
The coloring complex of $G$ is homotopy equivalent to a wedge of 
spheres, where the number of spheres is the number of acyclic 
orientations of $G$, and each sphere has dimension $n-3$.
\end{theorem}

\proof
First notice that 
the number of regions into which $A_G$ subdivides the sphere is the 
number of acyclic orientations of $G$, since points in the same region 
are all linear 
extensions of the associated acyclic orientation.  Therefore, $\Delta_G$ is the codimension one skeleton of a regular cell decomposition of an $(n-2)$-ball obtained by removing any single $(n-2)$-cell of $S^{n-2}.$  Since the ball has $A_G-1$ cells of dimension $n-2,$ its $(n-3)$-skeleton, and hence $\Delta_G,$ is homotopy equivalent to a wedge of $A_G-1$ spheres, all of dimension $n-3.$
\EOP

Jonsson also proved that $\Delta_G$ is constructible, and hence Cohen-Macaulay.  As we will see below, $\Delta_G$ has a convex ear decomposition which implies, by \cite[Theorem 4.1]{Sw}, that it is in fact doubly Cohen-Macaulay.  Specifically, if we remove any vertex from $A_G$ it remains an $(n-2)$-dimensional Cohen-Macaulay complex.

The arrangements viewpoint on the coloring complex follows easily from
a connection between bar resolutions and arrangements as developed in [HRW] and 
further exploited in [HW] and [PRW].  In particular, [HRW] deals with rings in which one
mods out by ideals in exactly the way the coloring complex arises, and [HRW] 
makes the connection in its more general setting to arrangements.

\section{Convex ear decomposition for the coloring complex}  \label{shell+ce}

The following notion was introduced by Chari in [Ch].

\begin{definition} Let $\Delta$ be a $(d-1)$-dimensional simplicial complex. 
A {\it convex ear decomposition} of  $\Delta $ is an ordered 
sequence $\Delta_1,\dots ,\Delta_m$ of pure $(d-1)$-dimensional subcomplexes
of $\Delta $ such that 
\begin{enumerate}
\item
$\Delta_1$ is the boundary complex of a $d$-polytope.  For each $j\ge 2$, 
$\Delta_j$ is a $(d-1)$-ball which is a proper subcomplex of the boundary of a 
simplicial $d$-polytope.
\item 
For $j\ge 2$, $\Delta_j \cap (\cup_{i<j} \Delta_i ) = \partial \Delta_j$.
\item
$\displaystyle\bigcup_j \Delta_j = \Delta $.
\end{enumerate}
\end{definition}
The subcomplexes $\Delta_1, \dots, \Delta_m$ are the {\it ears} of the decomposition.
The key ingredient in proving our main result is the lemma stated next, after
requisite terminology is introduced.  An 
arrangement $A = \{H_1,\dots ,H_s\} $ is {\it central} 
if each $H_i$ includes the origin, and $A$ is {\it essential } if 
$\cap_{i=1}^s H_i $ consists of exactly one point.
  For $A$ any essential central arrangement in $\real^n$, 
a {\it polytopal 
realization} of $A\cap S^{n-1}$ is any $n$-polytope containing the origin whose face fan is the fan of the arrangement.  Polytopal realizations of $A$ can be constructed by taking the polar dual of  Minkowski sums of line segments through the origin perpendicular to the hyperplanes (see, for instance, [Zi]).        

\begin{lemma} [Sw, Lemma 4.6]  \label{lemma}
         Let $A = \{H_1, \dots, H_s\}$ be an essential arrangement of hyperplanes in $\real^n.$  Let $P$ be any $n$-polytope whose face fan is the fan of $A.$ Let $H^+_{i_1}, \dots, H^+_{i_t}$ be  closed half-spaces of distinct hyperplanes in $A.$ If $B = \partial P \cap H^+_{i_1} \cap \dots \cap H^+_{i_t}$ is nonempty, then $\partial B$ is combinatorially equivalent to the boundary of an $(n-1)$-polytope.      \end{lemma}

\begin{theorem}\label{conv-ear}
The coloring complex of a graph has a convex ear decomposition.   Moreover, any
simplicial complex 
obtained by replacing $A_{K_n} $ in Theorem ~\ref{arrangement} by an essential,
central, simplicial arrangement and $A_G$ by any subarrangement will have a 
convex ear decomposition.\end{theorem}

\proof
Suppose that $G$ is connected.  Then  $A_G \cap H$ is an essential arrangement. Let $P$ be a polytopal realization of $S^{n-2} \cap A_G,$ and let $F_1, F_2, \dots, F_t$ be a line shelling of the facets of $P$ (as in e.g. [Zi]).  Identify each  facet with the corresponding  region of $A_G \cap S^{n-2}$ and, after further subdivision, a subcomplex of $A_{K_n} \cap S^{n-2}.$  By the lemma (applied in $A_{K_n} \cap S^{n-2}$), the boundary of each such region   is combinatorially equivalent to the boundary of a simplicial polytope.  Theorem \ref{arrangement} and the properties of line shellings imply that setting $\Delta_1 = \partial F_1,$ and for $2 \le i \le t-1, \Delta_i$  equal to the closure of $ \partial F_i \setminus (\partial F_1 \cup \dots \cup \partial F_{i-1}),$  produces a convex ear decomposition of $\Delta_G.$ 

For general finite graphs $G$, the intersection of all of the hyperplanes in $A_G$ is a $k$-dimensional subspace of $\real^n$, where $k$ is the number of components of $G.$  The lemma still implies that as a subcomplex of $A_{K_n} \cap S^{n-2}$ the boundary of each region of $A_G \cap S^{n-2}$ is combinatorially equivalent to the boundary of a simplicial polytope.  Let $H^\prime$ be the subspace of $\real^n$ orthogonal to the intersection of all of the hyperplanes in $A_G.$  Then the collection $A^\prime = \{H_1 \cap H^\prime, H_2, \cap H^\prime, \dots, H_s \cap H^\prime\},$ where the $H_i$ are the hyperplanes in $A_G,$ is an essential arrangement in $H^\prime.$  The facets of a polytopal realization of $A^\prime$ correspond to the regions of $A_G \cap S^{n-2}.$  Order the regions of $A_G \cap S^{n-2}$ in a way which corresponds to a line shelling of a polytopal representation of $A^\prime.$  Proceeding as before gives a convex ear decomposition of $\Delta_G.$  Indeed, the ea
 rs (and their intersections) are $(k-1)$-fold suspensions of a convex ear decomposition of the codimension one skeleton of a polytopal representation of $A^\prime.$  
 
 The only property of $A_{K_n}$ used above was the fact that it was a simplicial arrangement,
 so the above proof carries over immediately to the more general setting.  
\EOP

\begin{remark}
When $G$ is connected, the above reasoning also leads to an obvious shelling of $\Delta_G.$  However, the question of shellability is  more subtle for $G$ having $k>1$ components 
since not all the facets of the coloring complex actually intersect with the perpendicular 
space $H'$ to the $k$-dimensional space $U$ shared by all the hyperplanes in $A_G$.  See [Hu] for a shelling of the coloring complex for any $G$.
\end{remark}

\section{The unipolar  complex of a graph}

The unipolar complex of $G$ was introduced by Jonsson in \cite{Jo}.  Let $v_i$ be a vertex of $G.$  The {\it unipolar complex} of $G$ at $v_i$, denoted $ \Delta_{G(v_i)},$ is defined to be the subcomplex of $G_\Delta$ consisting of faces $\sigma$ such that $v_i \notin \bigcup^{r-1}_{j=1} S_j,$ where $S_1|\dots|S_{r-1}$ is the  ordered partition associated to $\sigma$.  From the arrangements point of view, $\Delta_{G(v_i)}$ may be realized  by taking the restriction of $\Delta_G$ to the intersection of half spaces of
the form $x_j \le x_i$ for all $j \neq i$.  It is easy to see that
this is still a simplicial complex and is  the codimension 
one skeleton of a pure subcomplex of the boundary of a convex polytope.

Jonsson proved that $\Delta_{G(v_i)}$ is constructible, hence  Cohen-Macaulay.  In general, it does not have a convex ear decomposition.  For instance, if $G$ is not connected, then any unipolar complex of $G$ is contractible, which is impossible for complexes with a convex ear decomposition.  However, if $v_i$ has degree $n-1,$ then we have the following.

\begin{theorem} \label{unipolar}
Let $v_i$ be a vertex of degree $n-1$ in $G.$  Then the unipolar complex of $G$ at $v_i$ has a convex ear decomposition.
\end{theorem}

\proof
As noted above, $\Delta_{G(v_i)}$ is the restriction to $A_G$ of the codimension one skeleton of the subcomplex of $A_{K_n}$ given by restriction to the half-planes $x_i \ge x_j.$  Since $v_i$ is incident to every vertex of $G,$ this is actually a subdivision of a subcomplex of $A_G.$  The proof of the lemma (see \cite{Sw}) shows that there is a point in $\real^n$ which ``sees'' only the regions of the aforementioned subcomplex of $A_G.$  Hence, there is a line shelling of a polytopal realization of $A_G$ such that the regions of the subcomplex are first.  Now we can use exactly the same reasoning as in the connected case of Theorem \ref{conv-ear}.
\EOP

\begin{remark}
When $v_i$ has degree $n-1,$ the above reasoning leads to an obvious shelling of $\Delta_{G(v_i)}.$
\end{remark}

\begin{question}
For which pairs $(G,v_i)$ does $\Delta_{G(v_i)}$ have a convex ear decomposition?
\end{question}

\section{Enumerative consequences}\label{enum-consequences}

The following connection between the coloring complex $\Delta_G$ and the chromatic polynomial
$P_G(t)$ was first given in [Ste].

\begin{theorem} \label{color h-vector} \cite{Ste}
Let $\Delta_G$ be the coloring complex of $G$ and let the $h$-vector of the double cone of $\Delta_G$ be $(h_0, \dots, h_n).$ Then
\begin{equation} 
(1-t)^n \sum^\infty_{j=0} [(j+1)^n - P_G(j+1)] t^j = h_0 + h_1 t + \dots + h_n t^n.
\end{equation}
\end{theorem}

Similarly, the $h$-vector of a unipolar complex can be computed from $P_G.$  Interestingly, it does not depend on the choice of vertex.

\begin{theorem} \label{unipolar h-vector} \cite[Theorem 2.5]{Jo}
Let $\Delta_G$ be the coloring complex of $G$ and let $(h^\prime_0, \dots, h^\prime_{n-2})$ be the $h$-vector of a unipolar complex of $\Delta_G$.  Then
\begin{equation} 
(1-t)^{n-1} \sum^\infty_{j=0} \frac{(j+1)^n - P_G(j+1)}{j+1} t^j = h^\prime_0 + h^\prime_1 t + \dots + h_{n-2} t^{n-2}.
\end{equation}
\end{theorem}

Since the $h$-vector of a cone equals the $h$-vector of the original complex, $h_{n-1}=h_n=0.$  In
order to state the enumerative consequences of Theorems \ref{conv-ear} and \ref{unipolar}, 
we first recall the definition of an M-vector.  

\begin{definition}  A sequence of  nonnegative integers $(h_0,h_1, \dots, h_d)$ is an {\bf M-vector} if it is the Hilbert function of a homogeneous quotient of a polynomial ring.  Equivalently, the terms form a degree sequence of an order ideal of monomials.
\end{definition}

Another definition given by arithmetic conditions is due to Macaulay.  Given  positive integers $h$ and $i$
there is a  unique way of writing
$$h=\binom{a_i}{i} + \binom{a_{i-1}}{i-1} + \dots + \binom{a_j}{j}$$
so that $a_i > a_{i-1} > \dots > a_j \ge j \ge 1.$  Define
$$h^{<i>} = \binom{a_i+1}{i+1} + \binom{a_{i-1}+1}{i} + \dots + 
\binom{a_j+1}{j+1} .$$

\begin{theorem} \cite[Theorem 2.2]{Sta} \label{CM h-vectors}
  A sequence of nonnegative integers $(h_0, \dots, h_d)$ is an M-vector if and only if $h_0=1$ and 
  $h_{i+1} \le h^{<i>}_i$ for all $1 \le i \le d-1.$
\end{theorem}

\begin{theorem}  \label{convear h-vector}
Suppose $\Delta$ is a $(d-1)$-dimensional complex with a convex ear decomposition.  Then,

\begin{enumerate}
  \item
    $h_0 \le h_1 \le \dots \le h_{\lfloor d/2 \rfloor}.$
  \item
    For $i \le d/2, \ h_i \le h_{d-i}.$
  \item
    $(h_0, h_1 - h_0, \dots, h_{\lceil d/2 \rceil} - h_{\lceil d/2 \rceil -1})$ is an M-vector.
 \end{enumerate}
 \end{theorem}

\proof  The first two inequalities are due to Chari \cite{Ch}.  The last statement is in \cite{Sw}.
\EOP

\begin{theorem} \label{main-color}
Let $G$ be a graph with $n$ vertices.  Define $h_0, \dots, h_n$ by the generating function equation
$$h_0 + h_1 t + \dots + h_n t^n = (1-t)^n \sum^\infty_{j=0} [(j+1)^n - P_G(j+1)] t^j.$$
Then
\begin{enumerate}
  \item  \label{main1}
    $h_0 \le h_1 \le \dots \le h_{\lfloor (n-2)/2 \rfloor}.$
  \item  \label{main2}
     For $i \le (n-2)/2, \ h_i \le h_{n-2-i}.$
   \item  \label{main3}
   $(h_0, h_1 - h_0, \dots, h_{\lceil (n-2)/2 \rceil} - h_{\lceil (n-2)/2 \rceil -1})$ is an M-vector.
\end{enumerate}
\end{theorem}

\proof
Theorems \ref{conv-ear}, \ref{color h-vector} and \ref{convear h-vector}.
\EOP

\begin{theorem}  \label{main-unipolar}
Let $G$ be a graph with $n$ vertices.  Suppose $G$ is chromatically equivalent to a graph which contains a vertex of degree $n-1.$  Define $(h^\prime_0, \dots, h^\prime_n)$ by the generating function formula
$$h^\prime_0 + h^\prime_1 t + \dots + h^\prime_{n-2} t^n = (1-t)^{n-1} \sum^\infty_{j=0} \frac{(j+1)^n - P_G(j+1)}{j+1} t^j.$$
Then
\begin{enumerate}  
  \item 
    $h^\prime_0 \le h^\prime_1 \le \dots \le h^\prime_{\lfloor (n-2)/2 \rfloor}.$
  \item  
     For $i \le (n-2)/2, \ h^\prime_i \le h^\prime_{n-2-i}.$
   \item 
   $(h^\prime_0, h^\prime_1 - h^\prime_0, \dots, h^\prime_{\lceil (n-2)/2 \rceil} - h^\prime_{\lceil (n-2)/2 \rceil -1})$ is an M-vector.
\end{enumerate}
\end{theorem}

\proof
Theorems \ref{unipolar}, \ref{unipolar h-vector} and \ref{convear h-vector}
\EOP

Let $A$ be a subarrangement of the $\mathcal{B}_n$ arrangement. The $\mathcal{B}_n$ arrangement consists of all the hyperplanes in $A_{K_n}$ and all coordinate hyperplanes $x_i=0.$ In \cite{Hu} Hultman proved the following relationship between $\chi_A(t),$ the characteristic polynomial of $A$ viewed as a matroid, and $(h^{\prime\prime}_0, \dots, h^{\prime\prime}_{n-1}),$ the $h$-vector of $\mathcal{B}_n \cap S^{n-1}$ restricted to $A.$

\begin{theorem} \cite{Hu} 
Let $A$ be a subarrangement of $\mathcal{B}_n$ and let $r$ be the rank of $A$ as a matroid.  Then
\begin{equation} \label{Bn-h}
  h^{\prime\prime}_0 + \dots + h^{\prime\prime}_{n-1} t^{n-1} = (1-t)^n \displaystyle\sum^\infty_{j=0} [(2j+1)^n - \chi_A(2j+1) (2j+1)^{n-r}] t^j.
\end{equation}

\end{theorem}

Combining Theorem \ref{conv-ear}, Theorem \ref{convear h-vector} and (\ref{Bn-h}) we obtain the following.

\begin{theorem} \label{main Bn}
Let $A$ be a subarrangement of $\mathcal{B}_n.$ Define $(h^{\prime\prime}_0, \dots, h^{\prime\prime}_{n-1})$ by (\ref{Bn-h}).  Then

\begin{enumerate}  
  \item 
    $h^{\prime\prime}_0 \le h^{\prime\prime}_1 \le \dots \le h^{\prime\prime}_{\lfloor (n-1)/2 \rfloor}.$
  \item  
     For $i \le (n-1)/2, \ h^{\prime\prime}_i \le h^{\prime\prime}_{n-1-i}.$
   \item 
   $(h^{\prime\prime}_0, h^{\prime\prime}_1 - h^{\prime\prime}_0, \dots, h^{\prime\prime}_{\lceil (n-1)/2 \rceil} - h^{\prime\prime}_{\lceil (n-1)/2 \rceil -1})$ is an M-vector.
\end{enumerate}

\end{theorem}

\begin{remark}
Characteristic polynomials of subarrangements of $\mathcal{B}_n$ correspond to chromatic polynomials of signed colorings introduced by Zaslavsky.  See \cite{Za}.

\end{remark}

In order to apply these methods to other arrangements it is essential that  subarrangements with the same characteristic polynomial (as matroids) have the same $h$-vector when restricted to the unit sphere.  In particular, all the simplicial subdivisions of the codimension one spheres corresponding to the hyperplanes must have the same $h$-vector.  

\begin{question}
Are there other (classes of) hyperplane arrangements such that the $h$-vectors of subcomplexes induced by subarrangements only depend on the characteristic polynomials of the subarrangements?  

\end{question}

\section{Matroids} \label{matroids}

Given the close connection between the chromatic polynomial of
a graph and the characteristic polynomial of the associated 
cycle matroid, it does not seem unreasonable to hope that it is possible to generalize Theorem \ref{main-color}  or Theorem \ref{main Bn} to matroids.  However, as the examples below show, it is not clear that there is 
any large class of matroids for which this is possible, though it is 
certainly possible that there is.

In these examples we let $\chi_M(t)$ be the characteristic polynomial of the matroid $M.$  When $G$ is connected, $P_G(t) = t \chi_{M_G}(t),$ where $M_G$ is the cycle matroid of the graph.  We will therefore use
\begin{equation} \label{matroid h}
h_0 + h_1 t + \dots + h_n t^n = (1-t)^n \sum^\infty_{j=0} [(j+1)^n - (j+1) \chi_M(j+1)] t^j.
\end{equation}
as the analog of the $h$-vector of the coloring complex for a rank $n-1$ matroid $M.$  

Let us now give examples violating various parts of Theorem \ref{main-color}.

\begin{example}
Let $M$ be $PG(5,2),$ the matroid whose elements correspond to the nonzero elements of the five-dimensional vector space over the field of cardinality two with their natural independence relations.  Then $\chi_M(t) = t^5 - 31 t^4 + 310 t^3-1240 t^2+1984 t-1024.$  Like the matroid associated to the braid arrangements, $M$ is binary and supersolvable.  However, (\ref{matroid h}) gives, $h_3 = -1678,$ a negative integer.

\end{example}

\begin{example}
Let $M$ be the matroid associated to the $B_3$ arrangement, the hyperplanes fixed by the symmetries of the cube.  Like the braid arrangements, $B_3$ is a free arrangement associated to a root system.   $\chi_M(t) = t^3 - 9 t^2 + 23 t - 15.$  Using (\ref{matroid h}) we find that $h_0 = 1, h_1 = 6, h_2= 47.$  The $h_i$ are nonnegative, but do not form an M-vector. 
\end{example}

\begin{example}
Let $\chi_M(t) = (t-1)^3 (t-2)(t-8)(t-10).$  Then $\chi_M(t)$ is the characteristic polynomial of the direct sum of $2$ coloops and the parallel connection of a $3$-point line, $9$-point line, and an $11$-point line \cite[Cor. 4.7]{Br}.  Now we find 
$$(h_0,\dots,h_5) = (1,121,472,4424,9167,2375).$$

\noindent This is an M-vector and satisfies (\ref{main1}) and (\ref{main2}) of Theorem \ref{main-color}.  However, (\ref{main3}) is not satisfied as 

$$(1,120,351,3952)$$ 
is not an M-vector.

\end{example}

Since every $A_G$ is a subarrangement of the $\mathcal{B}_n$ arrangement, characteristic polynomials of graphic matroids must satisfy Theorem \ref{main Bn}.  Perhaps this possibly weaker condition is satisfied by all matroids.  However, this is also not true.  

\begin{example}
Let $M$ be the matroid of PG(2,6).  Using (\ref{Bn-h}) as a definition with $n=6,$ we obtain $h_1 = -3047$ and $h_3 = -65638.$ 
\end{example}

Let us conclude by mentioning one class of matroids closely related to 
graphic matroids to which Theorem \ref{main-color} or \ref{main Bn} could 
perhaps apply.

\begin{question}
Let $M$ be a  regular matroid, namely a matroid representable over
every field.  Does $M$ satisfy either 
Theorem \ref{main-color} or Theorem \ref{main Bn}?
\end{question}

\end{document}